\documentclass{singlecol-new}

\usepackage{amsmath,amsfonts,amssymb}
\usepackage{url}
\usepackage{graphicx}

\theoremstyle{TH}{
\newtheorem{theorem}{Theorem}[section]

}

\theoremstyle{THrm}{
\newtheorem{remark}[theorem]{Remark}

}

\begin{document}

\thispagestyle{empty}

\setcounter{page}{1}

\LRH{P.A.F. Cruz, D.F.M. Torres and A.S.I. Zinober}

\RRH{A non-classical class of variational problems}

\subtitle{}

\title{A non-classical class of variational problems}

\authorA{Pedro A. F. Cruz}

\affA{Department of Mathematics,\\ 
University of Aveiro,\\ 
3810-193 Aveiro, Portugal\\ 
E-mail: pedrocruz@ua.pt}

\authorB{Delfim F. M. Torres*}

\affB{Department of Mathematics,\\ 
University of Aveiro,\\
3810-193 Aveiro, Portugal\\
E-mail: delfim@ua.pt\\ 
{*}Corresponding author}

\authorC{Alan S. I. Zinober}

\affC{Department of Applied Mathematics,\\
The University of Sheffield,\\
Sheffield, S10 2TN, UK\\
E-mail: a.zinober@sheffield.ac.uk}

%--------------------------------------------------

\begin{abstract} 
We study a new non-classical class of variational problems 
that is motivated by some recent research on the non-linear 
revenue problem in the field of economics. This class of 
problem can be set up as a maximising problem 
in the Calculus of Variations (CoV) or Optimal Control.
However, the state value at the final fixed time, $y(T)$, is {\em a priori\/}
unknown and the integrand is a function of the unknown $y(T)$.
This is a non-standard CoV problem.
In this paper we apply the new costate boundary conditions $p(T)$
in the formulation of the CoV problem. We solve some sample examples in this problem class
using the numerical shooting method to solve the resulting TPBVP, and incorporate the free $y(T)$
as an additional unknown.
Essentially the same results are obtained using symbolic algebra software.
\end{abstract}

%-----------------------------------------

%\noindent \textbf{2010 Mathematics Subject Classification}: 49K05, 49M05.

%-----------------------------------------

\KEYWORD{non-linear revenue problems;
transversality conditions; shooting method;
calculus of variations; optimal control.}

\REF{to this paper should be made as follows: 
Cruz, P.A.F., Torres, D.F.M. and Zinober, A.S.I. (2010) 
`A non-classical class of variational problems', 
{\it Int. J. Mathematical Modelling and Numerical Optimisation}, Vol.~1, No.~1, 
pp.xxx\textendash xxx.}

\begin{bio}
Pedro A. F. Cruz is Assistant Professor in the Department
of Mathematics at the University of Aveiro.
He obtained a Computational Science MSc from Instituto Superior T\'{e}cnico
and a PhD in Mathematics from the University of Aveiro.
He has published journal and conference works
in statistics, optimization, and control theory.
His main research interests are in the field of computational optimization.\vs{9}

\noindent Delfim F. M. Torres is Associate Professor of Mathematics
at the University of Aveiro; Scientific Coordinator of the
Control Theory Group (cotg); Editor-in-Chief of IJMS
and IJAMAS. He received the \emph{Licenciatura} degree
from the University of Coimbra,
and the MSc and PhD degrees from the University of Aveiro.
Delfim F. M. Torres has written more than 150 scientific and
pedagogical publications, and held positions 
of Invited and Visiting Professor
in several countries in Europe, Africa, Caucasus, and America.
His research interests include several topics 
in the areas of the Calculus of Variations and Optimal Control.\vs{9}

\noindent Alan Zinober is Professor of Nonlinear Control Theory
in the Department of Applied Mathematics at the University of Sheffield.
After obtaining the BSc and MSc degrees at the University of Cape Town
he was awarded the PhD degree from the University of Cambridge.
He has been the recipient of a number of Engineering
and Physical Sciences Research Council and other research grants.
He has published many journal and conference publications,
and has edited three research monographs.
The central theme of his research is in the field of
sliding mode control and other areas of nonlinear control theory.
\end{bio}

%--------------------------------------------------

\maketitle

%--------------------------------------------------

\section{Introduction}

The first result of the calculus of variations
ever discovered must have been the statement that the shortest
path joining two points is a straight line segment. Another
classical variational problem consists in finding,
amongst all simple closed plane curves of a given fixed
length, one that encloses the largest possible area.
It is well known since ancient times that the circle
is the shape that encloses maximum area for a given length of
perimeter. However, it was not until the eighteenth century that a
systematic theory, the Calculus of Variations (CoV), began to emerge.
A modern face to the CoV is given by the theory of optimal control.
Economics is a source of interesting applications
of the theory of calculus of variations and optimal control.
Classical examples include the optimal capital spending problem,
optimal reservoir control, optimal production subject to
royalty payment obligations, optimal maintenance and replacement policy,
and optimal drug bust strategy \cite{Ngo,Sethi}.

The following economics problem (explained briefly here) has motivated this paper \cite{Zinober}.
A standard feature of the theory of the firm is that a profit
maximising firm facing a downward sloping demand curve
reacts to an increase in marginal cost by reducing output
and increasing price. In this context, it is well understood
that a requirement to pay a flat-rate royalty on sales
has just this effect of increasing
marginal cost and thereby decreasing output while
simultaneously increasing price.
However, the effect of permitting the royalty to take on more
general forms leads naturally to non-standard CoV problems,
and explains why this question has remained unaddressed to date
\cite{Zinober}. Recently the effect of piecewise linear
cumulative royalty schedules on the optimal intertemporal
production policy, \textrm{i.e.}, an optimal economics control problem
that does not fit into the classical class of variational problems,
has been formulated \cite{Zinober}.
The economics problem lies in the area of Repayable Launch Investment (RLI).
For the purposes of this paper we will outline just the mathematical nature
of the problem since the precise (nonlinear) economic details are of secondary importance here.
Consider the system in the time domain modelled by the differential equation
$$ y'(t) = u(t)\, , \qquad y(0) \mbox{ \ known}$$
with the endpoint state value $y(T)$ at time $t = T$ unknown.
We wish to determine the control function $u(t)$ for $ t \in [0, T]$
that maximises the return
$$
J[u(\cdot)] = \int_0^T f(t,y(t),u(t),y(T)) dt \, .
$$
Note that the integrand depends upon the {\it a priori\/} unknown final value $y(T)$.
This class of problem is not contained within the classical class of variational problems and the present paper indicates how such problems may be solved.

The manuscript is organized as follows. In Section~\ref{sec:nec:cond}
we develop the necessary conditions for the extremising solution.
The main idea is borrowed from Malinowska and Torres \cite{Malinowska:Torres},
where appropriate natural boundary conditions are proved
for problems of the calculus of variations on time scales
\cite{Ferreira:Torres}. We then consider two approaches to obtain the solution
of a sample example with a continuous integrand $f$ (Section~\ref{sec:ex}).
The first approach considers the numerical shooting solution (Section~\ref{sec:NA}).
The results obtained are then validated by symbolic algebra computations (Section~\ref{sec:SAS}).
We finish with conclusions in Section~\ref{sec:conc}.

%--------------------------------------------------

\section{The non-classical variational problem}
\label{sec:nec:cond}

We begin by developing the necessary conditions for
the extremising solution. Let $J$ be a functional of the form
$$
J[y(\cdot)] = \int_a^T f(t,y(t),y'(t),y(T)) dt
$$
where $$(t,y,y',z) \rightarrow f(t,y,y',z)$$ is
a smooth function and $T > a$. We consider the problem
of determining the functions $y(\cdot) \in C^1$ such that
$J[\cdot]$ has an extremum. An initial condition
$y(a) = \alpha$
is imposed on $y(\cdot)$, but $y(T)$ is free.

Suppose that $J[\cdot]$ has an extremum at $\tilde{y}(\cdot)$.
We can proceed as Lagrange did (\textrm{cf.} \cite{Gelfand}),
by considering the value of $J$ at a nearby function
$y = \tilde{y} + \varepsilon h$, where
$\varepsilon$ is a small parameter, $h(\cdot) \in C^1$, and $h(a) = 0$.
Because $y(T)$ is free, we do not require $h(\cdot)$ to vanish at $T$. Let
$$
\phi(\varepsilon) = J[(\tilde{y} + \varepsilon h)(\cdot)]
= \int_a^T f(t,\tilde{y}(t)+\varepsilon h(t),\tilde{y}'(t)
+ \varepsilon h'(t),\tilde{y}(T)+\varepsilon h(T)) dt \, .
$$
A necessary condition for $\tilde{y}(\cdot)$ to be an
extremizer is given by
\begin{equation}\label{eq:FT}
\left.\phi'(\varepsilon)\right|_{\varepsilon=0} = 0
\Leftrightarrow \int_a^T \left[
f_y(\cdots) h(t) + f_{y'}(\cdots) h'(t)
+ f_z(\cdots) h(T) \right] dt = 0 \, ,
\end{equation}
where $(\cdots) = \left(t,\tilde{y}(t),\tilde{y}'(t),\tilde{y}(T)\right)$.
Integration by parts gives
$$
\int_a^T f_{y'}(\cdots) h'(t) dt
= \left. f_{y'}(\cdots) h(t)\right]_a^T
- \int_a^T \frac{d}{dt} \left(f_{y'}(\cdots)\right) h(t) dt \, .
$$
Because $h(a) = 0$, the necessary condition (\ref{eq:FT}) can
be then written as
\begin{multline}
\label{eq:aft:IP}
0 = \int_a^T \Biggl\{ \left[f_y(\cdots)-\frac{d}{dt}
f_{y'}(\cdots)\right] h(t) \\
+ \left[\frac{f_{y'}\left(T,\tilde{y}(T),\tilde{y}'(T),
\tilde{y}(T)\right)}{T - a}
+ f_z(\cdots)\right] h(T) \Biggr\} dt
\end{multline}
for all $h(\cdot) \in C^1$ such that $h(a) = 0$.
In particular, equation (\ref{eq:aft:IP}) holds for the
subclass of functions $h(\cdot) \in C^1$ that do vanish
at $h(T)$. Thus, the classical arguments apply, and therefore
\begin{equation}
\label{eq:EL}
f_y(\cdots)-\frac{d}{dt}
f_{y'}(\cdots) = 0 \, .
\end{equation}
Equation (\ref{eq:aft:IP}) must be satisfied for all $h(\cdot) \in C^1$
with $h(a) = 0$, which includes functions $h(\cdot)$ that do not vanish at $T$.
Consequently, equations (\ref{eq:aft:IP}) and (\ref{eq:EL}) imply that
\begin{gather*}
\int_a^T \left[\frac{f_{y'}\left(T,\tilde{y}(T),\tilde{y}'(T),
\tilde{y}(T)\right)}{T - a}
+ f_z(\cdots)\right] h(T) dt = 0\\
\Leftrightarrow
h(T) \left( f_{y'}\left(T,\tilde{y}(T),\tilde{y}'(T),
\tilde{y}(T)\right)
+ \int_a^T f_z(\cdots) dt
\right) = 0 \, ,
\end{gather*}
that is,
\begin{equation}
\label{new:nat:bc}
f_{y'}\left(T,\tilde{y}(T),\tilde{y}'(T),
\tilde{y}(T)\right)
+ \int_a^T f_z(\cdots) dt = 0 \, .
\end{equation}
We remark that in the classical setting $f$
does not depend on $y(T)$, that is,
$f_z = 0$. In that case (\ref{new:nat:bc})
reduces to the well known natural boundary
condition $f_{y'}\left(T,\tilde{y}(T),\tilde{y}'(T)\right) = 0$
(or, from an Hamiltonian optimal control perspective, $p(T) = 0$).
We have just proved the following result:

\begin{theorem}
\label{thm:mr}
Let $a$ and $T$ be given real numbers, $a < T$.
If $\tilde{y}(\cdot)$ is a solution of the problem
\begin{equation}
\label{eq:P}
\begin{gathered}
J[y(\cdot)] = \int_a^T f(t,y(t),y'(t),y(T)) dt \longrightarrow \textrm{extr} \\
y(a) = \alpha \quad (y(T) \text{ free}) \\
y(\cdot) \in C^1 \, ,
\end{gathered}
\end{equation}
then
\begin{equation}\label{eqn1}
\frac{d}{dt} f_{y'}\left(t,\tilde{y}(t),\tilde{y}'(t),\tilde{y}(T)\right)
= f_y\left(t,\tilde{y}(t),\tilde{y}'(t),\tilde{y}(T)\right)
\end{equation}
for all $t \in [a,T]$.
Moreover,
\begin{equation}\label{eqn2}
f_{y'}\left(T,\tilde{y}(T),\tilde{y}'(T),\tilde{y}(T)\right)
= - \int_a^T f_z\left(t,\tilde{y}(t),\tilde{y}'(t),\tilde{y}(T)\right) dt \, .
\end{equation}
\end{theorem}

\begin{remark}
From an optimal control perspective one has
\begin{equation*}
p(T) = f_{y'}\left(T,\tilde{y}(T),\tilde{y}'(T),
\tilde{y}(T)\right) \,,
\end{equation*}
where $p(t)$ is the Hamiltonian multiplier.
Theorem~\ref{thm:mr} asserts
that the usual necessary optimality conditions
(the Euler-Lagrange equation \cite{Gelfand} or the
Pontryagin maximum principle \cite{Pontryagin})
hold for problem (\ref{eq:P}) by substituting
the classical transversality condition $p(T) = 0$ with
\begin{equation}
\label{AZ}
p(T) = - \int_a^T f_z\left(t,\tilde{y}(t),\tilde{y}'(t),\tilde{y}(T)\right) dt \, .
\end{equation}
\end{remark}

% ------------------------------------------

\section{An illustrative example}
\label{sec:ex}
We consider an example that illustrates the new class of CoV problem.
It has the same form as the complicated nonlinear optimal intertemporal
production policy problem.
Consider the ODE system described by
\begin{equation}
y'(t) = u(t)\, , \qquad y(0)=0\, . \label{A1}
\end{equation}
We wish to maximise
\begin{equation}\label{eq:functional}
J[u(\cdot)] = \int_0^T f(t,y(t),u(t),z)\, dt
\end{equation}
where
\begin{equation}\label{eq:lagrangean}
f(t,y,u,z) = a \sqrt{u}-\left(\frac{3}{4} + z\sin(\pi t/10)\right)u
\end{equation} is a continuous function.
The initial known state is $y(0)=0$ and final state value $ z = y(T) $ is free.
In this example we set $T=10$.
The Hamiltonian is $H(t,y,u,p) = -f + p \cdot u$ and
\begin{equation*}
\begin{cases}
y'(t) =    H_p(t,y(t),u(t),p(t)) \\
p'(t) =  - H_y(t,y(t),u(t),p(t))\,.
\end{cases}
\end{equation*}
Function $f$ does not depend on $y$, and for an optimum (maximum in this example),
the costate satisfies
\begin{equation}\label{A2}
p' = -H_y \Leftrightarrow  p' = 0 \, .
\end{equation}
The stationarity condition is
\[
  H_u = 0
\]
and this yields
\begin{equation} \label{A3}
u(t) = \frac{1}{4} (z \sin(\pi t/10) - p(t))^2 \,.
\end{equation}
From (\ref{AZ})
\[
p(T) = - \int_0^{10} f_z\left(t,\tilde{y}(t),\tilde{y}'(t),\tilde{y}(T)\right) dt
\]
holds, \textrm{i.e.},
\begin{equation}
\label{A4}
p(T) = \int_0^{10} \sin(\pi t/10) \, u(t) \, dt \, .
\end{equation}

% ------------------------------------------

\section{Numerical shooting algorithms}
\label{sec:NA}

Let us consider the necessary conditions (NC) that need to be satisfied.
For the system of ODEs (\ref{A1}) and (\ref{A2}) with control (\ref{A3}),
the known zero initial condition $y(0)$ and a guessed initial value $p(0)$,
we need to ensure that the natural boundary condition (\ref{A4}) is satisfied.

We need to solve the two point boundary value problem.
Also we need to iterate the value of $z$ used in (\ref{A3})
to ensure that in fact the value $z$ equals the value  obtained
for $y(t)$ at $t=T$. When one has obtained convergence regarding the values $y(T)$
used in (\ref{A3}) and $p(T)$ (\ref{A4}),
then NC is satisfied and we should have the optimal solution.

Use the Newton shooting method with two guessed values $v_1$ and $v_2$ \cite{Betts}.
We desire $v_1 = p(0)$,  and $v_2 = p(T)$ as specified by equation (\ref{A4}).
When the program obtains results with these two equations holding to a very high degree
of accuracy, the necessary conditions NC hold and we should have the optimal solution.
We have solved the shooting method problem using \textsf{C++} and the highly accurate
Numerical Recipes library routines \cite{Press}:

\bigskip

(i) We integrate the system $(y(t), p(t), g_p(t), J(t))$, \textrm{i.e.},
the system of ODEs (\ref{A1}) and (\ref{A2}), and
$$ g_p'  =  \sin(\pi t/10) u(t) \, , $$
$$ J' = g \, .$$
The results are
$y(T)=0.86928249597392515$,
$p(T) = - 0.46111638323272386$,
$g_p(T) = - 0.46111638323273074$, and
\begin{equation}
\label{AZJ}
J(T) = 1.85448307363352  \, .
\end{equation}
Perturbations of the optimal control $u(t)$ by increasing and decreasing the
value of $u(t)$ at a single time instant yield smaller $J(T)$ values.
See Figure~1 for results on state variable $y(t)$ and control variable $u(t)$.

A completely different approach, using a nonlinear programming
technique, was also used. This technique may be useful for the actual 
piecewise constant economics problem. We solved this problem using Euler and Runge-Kutta
discretisation, and an optimisation algorithm to solve for the unknown control variables
$u_k$ at each time instant $t_k$. We computed the
nonlinear programming problem using AMPL \cite{AMPL}
with the MINOS solver and NEOS \cite{NEOS}.
Using 40 time steps yields a good approximation very similar
to the optimal results obtained using the precise approach here described.

% --------------------------------------------------------

\section{Symbolic algebra solution}
\label{sec:SAS}

Consider the ODE system (\ref{A1}) and the associated
optimal control problem described
by (\ref{eq:functional}) and (\ref{eq:lagrangean}).
We set this as a minimization problem. From  (\ref{eq:functional}) define
\begin{equation}
\label{eq:J}
J_m[u(\cdot)] = \int_0^T g(t,y(t),u(t),z)\, dt
\end{equation}
where
\begin{equation}\label{eq:f}
g(t,y,u,z) = -f(t,y,u,z) = \left(\frac{3}{4} + z\sin(\pi t/10)\right)u - a \sqrt{u}
\end{equation}
with the final state value $ z = y(T) $ free, and $T=10$.
We now use the Euler-Lagrange equation (\ref{eqn1}) to find candidate solutions:
\begin{equation*}
\frac{d}{dt} g_{y'}\left(t,\tilde{y}(t),\tilde{y}'(t),\tilde{y}(T)\right)
= g_y\left(t,\tilde{y}(t),\tilde{y}'(t),\tilde{y}(T)\right)
\end{equation*}
for all $t \in [0,T]$. Set $u=y'$ and $z=y(T)$ so
\begin{eqnarray*}
\frac{d}{dt} g_{y'} (t,y,y', y(T) )
  & = & \frac{d}{dt} g_{u} (t,y,u,z) \\
  & = & \frac{d}{dt} \left( \frac{3}{4} - \frac{1}{2\sqrt{u}} + z \sin\left(\frac{\pi t}{10}\right) \right) \\
  & = & \frac{1}{10} \pi z \cos\left(\frac{\pi t}{10}\right)+\frac{u'}{4 u^{3/2}} \, .
\end{eqnarray*}
Since $g_y=0$, using (\ref{eqn1}) we can find $u$ by solving
\begin{equation*}
\frac{1}{10} \pi z \cos\left(\frac{\pi t}{10}\right)+\frac{u'(t)}{4 u(t)^{3/2}} = 0 \,.
\end{equation*}
Next result was obtained using \textsf{Maple}:
\begin{equation*}
u(t)=\frac{1}{\left(c+2 z \sin \left(\frac{\pi  t}{10}\right)\right)^2} \, .
\end{equation*}
We find easily $y$ from $y'=u$ using integration, \textrm{e.g.}, in \textsf{Maple} or \textsf{Mathematica}:
\begin{equation}
\label{eq:yt}
y(t)=
\frac{20 \left(-\tan^{-1}
   \left(2 z / D\right)
   g(t) c^2+ h(t)
   g(t) c^2-z D
   \left(-\cos
   \left(\frac{\pi  t}{10}\right) c+ g(t) \right)\right)}{c \pi
   D^{3}
   g(t)}
\end{equation}
where
\begin{eqnarray*}
D    & = & \sqrt{c^2-4 z^2} \, , \\
g(t) & = & \left(c+2 z \sin \left(\frac{\pi  t}{10}\right)\right) \, ,\\
h(t) & = & \tan ^{-1}
   \left(\frac{2 z+c \tan
   \left(\frac{\pi  t}{20}\right)}{\sqrt{c^2-4 z^2}}\right) \, .
\end{eqnarray*}

Some comments:

\begin{itemize}

\item The function $h$ is not defined for $t=10$. However one can define $h$ for $t=10$ as
\[
\lim_{t\to T} \tan^{-1} \left(\frac{2 z+c \tan
   \left(\frac{\pi  x}{20}\right)}{\sqrt{c^2-4 z^2}}\right) = \frac{\pi}{2} \, .
\]

\item Theorem~\ref{thm:mr} assumes $u \in C^0$.

\item From (\ref{eq:f}) $u(t)\ge0$. If $u=0$ we have $J=0$
and we see that this is not the best solution.
So $u>0$ and $z=y(T)>0$.

\item We must verify $c^2-4 z^2>0$ (see $D$ in (\ref{eq:yt}))
so two cases are to be investigated: $c>2z$ and $c<-2z$.

\end{itemize}

Recall (\ref{eqn2}):
\begin{equation}
\label{eq:n7}
g_{y'}\left(T,\tilde{y}(T),\tilde{y}'(T),
\tilde{y}(T)\right) = - \int_a^T g_z\left(t,\tilde{y}(t),\tilde{y}'(t),\tilde{y}(T)\right) dt \, .
\end{equation}
We have
\begin{equation*}
g_z(t,y,u,z) = u \sin \left(\frac{\pi  t}{10}\right)
\end{equation*}
and using $u(\cdot)$ we arrive to
\begin{equation*}
g_z(t,y,u,z) = \frac{\sin \left(\frac{\pi  t}{10}\right)}{\left(c+2 z \sin \left(\frac{\pi t}{10}\right)\right)^2} \, .
\end{equation*}
Integrating for the branch $c>2z$ we obtain
\begin{equation*}
\int_0^{10} g_z\left(t,y(t),y'(t),y(T)\right) dt
= \frac{20 \left(2 \tan^{-1}\left(\frac{2 z}{D}\right) z-\pi  z+D \right)}{\pi  D^3} \, .
\end{equation*}
The left-hand side of (\ref{eq:n7}) is
\begin{equation*}
g_{y'}\left(T,y(T),y'(T),y(T)\right) = \frac{1}{4} (3-2 c) \, .
\end{equation*}
Solving equation (\ref{eq:n7}) numerically we get $c = 2.42223$ and $z= 0.869282$.
The objective value is $J_m=-1.85448$ obtained using numerical integration over (\ref{eq:J}).
This compares favorably with the result (\ref{AZJ}).  
Note that similar calculations for the branch $c<-2z<0$
provide a worse solution with $J_m=+6.62857$
($c = -7.21816$ and $z = 3.14287$).

The results obtained here by Symbolic Algebra Computations (SAC)
are in accordance and validate the numerical shooting solution
obtained in Section~\ref{sec:NA}.
\begin{figure}\label{fig_solution}
\includegraphics[width=6cm]{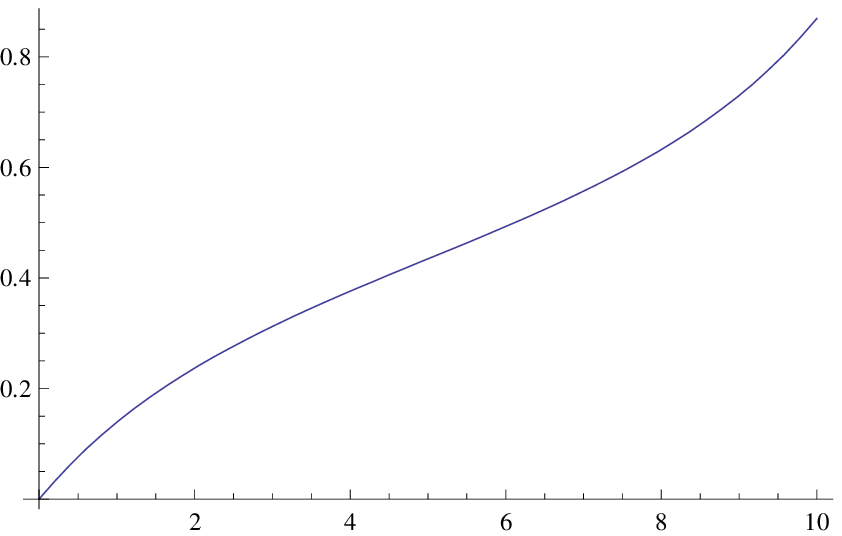}
\includegraphics[width=6cm]{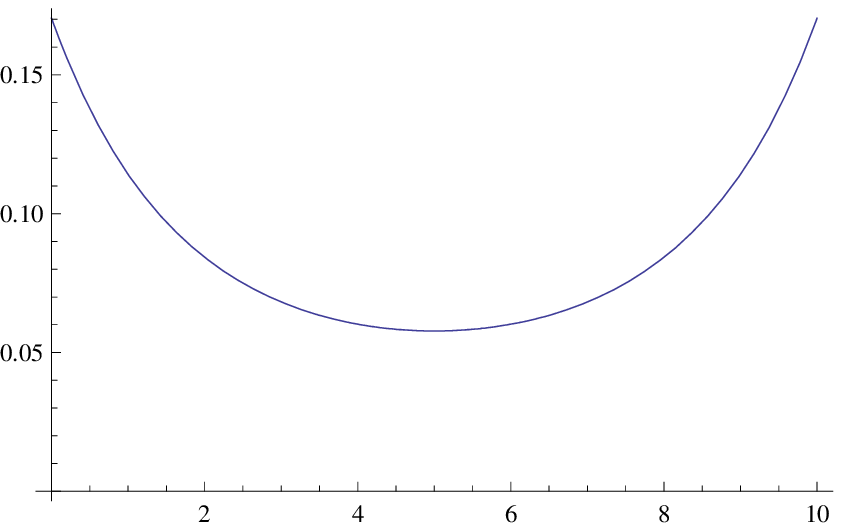}
\caption{Optimal pair $(y(t), u(t))$ to problem (\ref{A1})-(\ref{eq:functional}) 
obtained from both symbolic computation and the shooting method.}
\end{figure}

% ------------------------------------------

\section{Conclusion}
\label{sec:conc}

In this note we have shown how the standard necessary optimality conditions
and numerical procedures for problems of the calculus
of variations and optimal control
should be adapted in order to cover Lagrangians
depending on the free end-point. The numerical techniques
were validated with a simple sample example that allows
symbolic calculations using a modern computer algebra system.
In the actual optimal intertemporal
production policy economics problem the Lagrangian may be piecewise continuous
and this requires amended numerical techniques, such as nonlinear programming,
for its solution. This numerical solution approach will be important for solution 
of the actual nonlinear economics problem.

% ------------------------------------------

\section*{Acknowledgements}

We thank the reviewers for their helpful comments.
The first two authors were supported
by the Centre for Research on Optimization and Control (CEOC)
from the Portuguese Foundation for Science and Technology (FCT),
cofinanced by the European Community fund FEDER/POCI 2010.

%-----------------------------------------------

\NUMBIB

%-----------------------------------------------

\end{document}